\documentclass{svproc}

\PassOptionsToPackage{x11names}{xcolor}
\usepackage[dvipsnames]{xcolor}
\usepackage{graphicx}

\usepackage{amsmath,amsthm}
\usepackage{amssymb}
\usepackage{pgf}
\usepackage{tikz}
\usetikzlibrary{fadings, decorations.pathreplacing}\usetikzlibrary{arrows.meta}
\usepackage{multicol}
\usepackage{tcolorbox}
\usepackage{mathtools}
\let\llncssubparagraph\subparagraph
\let\subparagraph\paragraph
\usepackage{titlesec}
\let\subparagraph\llncssubparagraph
\usepackage{setspace}
\usepackage[T1]{fontenc}
\usepackage[percent]{overpic}
\usepackage{hyperref}
\usepackage{multirow}
\usepackage{tikz}

\usepackage{bbm}
\usepackage{cleveref}

\usetikzlibrary{arrows,automata}

\usepackage{tikz-network}

\usepackage[commandnameprefix=always]{changes}
\colorlet{change}{Firebrick3}
\colorlet{draft}{RoyalBlue3}

\usepackage{url}

\usepackage[utf8]{inputenc}
\usepackage{amsmath}
\usepackage{amsfonts}
\usepackage{amssymb}
\usepackage{amsthm}
\usepackage{graphicx}
\graphicspath{{./img}}
\usepackage{amsthm}
\usepackage{enumitem}
\theoremstyle{remark}

\usepackage{tcolorbox}
\newcommand{\code}[1]{\mbox{
    \ttfamily
    \tcbox[
        on line,
        boxsep=0pt, left=4pt, right=4pt, top=2pt, bottom=1.5pt,
        toprule=0pt, rightrule=0pt, bottomrule=0pt, leftrule=0pt,
        oversize=0pt, enlarge left by=0pt, enlarge right by=0pt,
        colframe=white, colback=black!12
    ]{#1}
}}

\usepackage{mathtools}
\usepackage{xcolor}
\usepackage{textpos}
\usepackage[percent]{overpic}

\usepackage[backend=bibtex,giveninits=true,style=numeric,citestyle=numeric, doi=false,isbn=false,url=false,eprint=false,maxbibnames=99]{biblatex}
\bibliography{main}

\begin{document}
\mainmatter              
\title{A Payne-Whitham model of urban traffic networks in the presence of traffic lights and its application to traffic optimisation}
\titlerunning{Payne-Whitham model of urban traffic in the presence of traffic lights}  
%
\author{Mauritz N. Cartier van Dissel\inst{1} \and Pawe\l\ Gora\inst{2} \and Dragos Manea\inst{3}}
\authorrunning{Cartier van Dissel, Gora and Manea} 
%
%
\institute{Complexity Science Hub, Vienna, Austria
\and
Faculty of Mathematics, Informatics and Mechanics, University of Warsaw, Warsaw, Poland
\and
``Simion Stoilom'' Institute of Mathematics of the Romanian Academy, Bucharest, Romania
}

\maketitle  

\begin{abstract}
    Urban road transport is a major civilisational and economic challenge, affecting the quality of life and economic activity. Addressing these challenges requires a multidisciplinary approach and sustainable urban planning strategies to mitigate the negative effects of traffic in cities. In this paper, we introduce an extension of one of the most popular macroscopic traffic simulation models, the Payne-Whitham model. We investigate how this model, originally designed to model highway traffic on straight road segments, can be adapted to more realistic conditions with arbitrary road network graphs and multiple intersections with traffic signals. Furthermore, we showcase the practical application of this extension in experiments aimed at optimising traffic signal settings. For computational reasons, these experiments involve the adoption of surrogate models for approximating our extended Payne-Whitham model, and subsequently, we utilise the Differential Evolution optimization algorithm, resulting in the identification of traffic signal settings that enhance the average speed of cars and decrease the total length of queues, thereby facilitating smoother traffic flow.
\end{abstract}

\section{Introduction}\label{sec:intro}
Urban road transport is a major civilisational and economic challenge, affecting quality of life and economic activity. Inefficiency in transport can lead to travel delays, driver stress, noise, wasted fuel and energy, air pollution, problems in organising public transport and detours, etc. Addressing these challenges requires a multidisciplinary approach and sustainable urban planning strategies to mitigate the negative effects of road transport in cities. Professional traffic simulators are among the standard tools used to plan and analyse transport in urban areas. They usually implement some transport demand models (specifying the number of trips between different urban areas), trip assignment models (assigning trips to specific routes), mode choice models (i.e., models of selecting means of transport), and finally, traffic models describing how traffic flows on the road network \cite{4stepmodel}. 

There are different approaches to traffic modelling and they can be classified according to different factors. One of the most popular classifications is based on the level of detail. There are macroscopic models which describe traffic from the perspective of the collective flow of vehicles, using aggregated values such as average speed, congestion or density. There are also microscopic models, which consider each vehicle individually and describe the time-space behaviour of individual drivers under the influence of vehicles in their vicinity. Finally, there are mesoscopic models, which describe traffic flow at an intermediate level of detail between microscopic and macroscopic models. In such models, vehicles and drivers' behaviours can be represented in groups, or they can still be described individually, but in a more aggregated way. Some studies also consider nanoscopic (or sub-microscopic) models, which contain mathematical representations of the vehicle subunits and their interactions (internally and with the environment).

Many traffic models have been developed over the years and a review of the most popular and successful models can be found in \cite{HOO}.
There is a famous quote attributed to George Box ``All models are wrong, but some are useful'', meaning that while models may not perfectly represent reality, they can still provide valuable insight and utility in understanding and predicting certain phenomena. The use of particular models usually depends on the purpose and scope of the modelling. For example, to design intersections or investigate road safety, it is usually necessary to model all interactions between vehicles, so the use of microscopic models seems appropriate. On the other hand, to estimate travel times, speeds or congestion, macroscopic or mesoscopic models may be sufficient.

In general, the advantage of microscopic models is a relatively good accuracy and the ability to model complex behaviours caused by interactions between vehicles. Examples are the Wiedemann model \cite{wiedemann1974simulation}, the Gipps model \cite{gipps1981behavioural} and the Nagel-Shreckenberg model \cite{nagel1992cellular}. Although these models give relatively good conformity with real-world conditions, they are often computationally expensive, making them sometimes difficult to use on a large scale (however, with increasing computing power and the ability to distribute computations, they become more accessible for large-scale traffic simulation, as demonstrated in the case of the Traffic Simulation Framework software \cite{gora2009traffic}). They are usually also difficult to calibrate, as they have more parameters (related to the behaviour of individual vehicles) than macroscopic or mesoscopic models, and often require large amounts of real traffic data. Therefore, it is also important to develop the other types of models.

Nanoscopic models can be even more computationally demanding, so they are usually used on a relatively small scale, such as a single road segment, a single intersection, or a small group of intersections. Examples are SIMONE \cite{minderhoud1999supported}, which is suitable for representing conventional motorway traffic and the foreseen future with vehicles equipped with Intelligent Cruise Control (ICC) systems, and MIXIC \cite{van1995}, which is also suitable for studying autonomous ICC systems.

Mesoscopic models, which are intermediate models between microscopic and macroscopic levels, are usually more computationally efficient and can also be applied to large-scale traffic simulations. Examples are the models developed by Prigogine and Herman \cite{prigogine1971kinetic, prigogine1961boltzmann}, which are based on the kinetic gas theory.

Macroscopic models are usually based on partial differential equations describing the evolution of aggregated values such as traffic density, traffic flow and average speed. Most of them are derived from the first and most popular macroscopic model, the Lighthill-Whitham model \cite{lighthill1955kinematic}. In this article, we examine the widely used macroscopic traffic simulation model, the Payne-Whitham model \cite{Payne1971,Whitham1999}, originally developed for simulating highway traffic on straight road segments. We extend this model to accommodate more realistic scenarios, incorporating arbitrary road network graphs and multiple intersections regulated by traffic signals. This enhanced version, referred to as the Payne-Whitham model with Traffic Lights (PWTL), is subsequently applied in experiments aimed at optimizing traffic signal configurations for improved traffic flow. 

The problem of traffic signal control has been extensively studied, with approaches ranging from static and actuated control of individual intersections using fixed rules (e.g., Webster's method \cite{webster}) to more advanced coordination over larger areas with multiple intersections. Advanced techniques include artificial intelligence (e.g., reinforcement learning \cite{el2013multiagent, wei2019colight} and evolutionary algorithms \cite{gora2019solving}) and quantum computing \cite{inoue2021traffic}. In our approach, we use the Differential Evolution algorithm \cite{storn1997differential}, which is well-suited to our traffic light optimization task due to its effectiveness in handling complex optimization problems. 

These optimization methods require a traffic model to evaluate the quality (measured by some traffic characteristics such as average speed, travel times, and delays) of traffic signal settings. Since such evaluations are usually time-consuming, surrogate models are often employed to approximate the results of the original models quickly and with good accuracy \cite{skowronek2021graph, skowronek2023graph}. The high computational demand was also an issue with the PWTL presented, so we tested four surrogate models to predict the total average speed and the overall length of the queues: Linear Regression, Multi-Layer Perceptron (MLP), Support Vector Regression (SVR), and Random Forest. Although Linear Regression was the least accurate surrogate model, it was still able to find traffic light configurations that optimized speed and queue length effectively.

The article is organised as follows: Section \ref{sec:model} summarises the original Payne-Whitham model, Section \ref{sec:finite} introduces a numerical approximation for the Payne-Whitham model on a single road and provides a stability condition for the time and speed discretisation parameters, Section \ref{sec:intersections} adds the treatment of intersections, Section \ref{sec:init} explains how we calibrated and initialised the model with real traffic data, Section \ref{sec:optim} introduces the traffic signal optimisation procedure, while the experiments with applications of the developed model to the traffic signal optimisation task are described in Section \ref{sec:optexp}. Finally, Section \ref{sec:conclusions} concludes the article.

\section{Macroscopic models -- Payne-Whitham model}\label{sec:model}
Developed in the 1970s \cite{Payne1971,Whitham1999}, the Payne-Whitham model is a macroscopic traffic simulation model that incorporates advanced features in terms of driver behaviour behind the wheel, which is influenced by the level of road congestion and other related factors. The model was first developed for a single road segment and later Kotsialos et. al. \cite{KotsialosPapageorgiou} came up with an implementation for road networks, paying special attention to the behaviour of cars at junctions. 

Being a macroscopic model, the Payne-Whitham model consists of an equation describing the joint evolution in time of the density (the number of cars per a small distance per lane) and average speed of cars at each point on the road. The fact that not only the density, but also the speed evolves according to a partial differential equation classifies this model as a \textit{second-order} traffic model (see \cite{CaligarisSaconeSiri2011} for details).

\begin{figure}
\center{\includegraphics[scale=0.37]{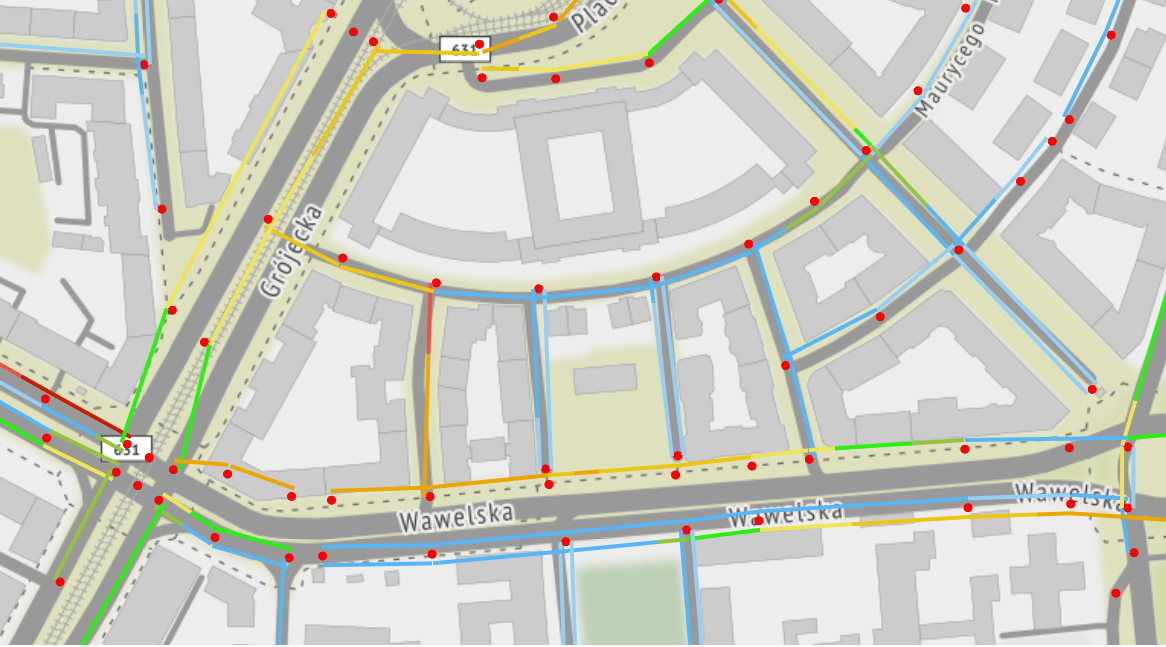}}
\center{\includegraphics[scale=0.06]{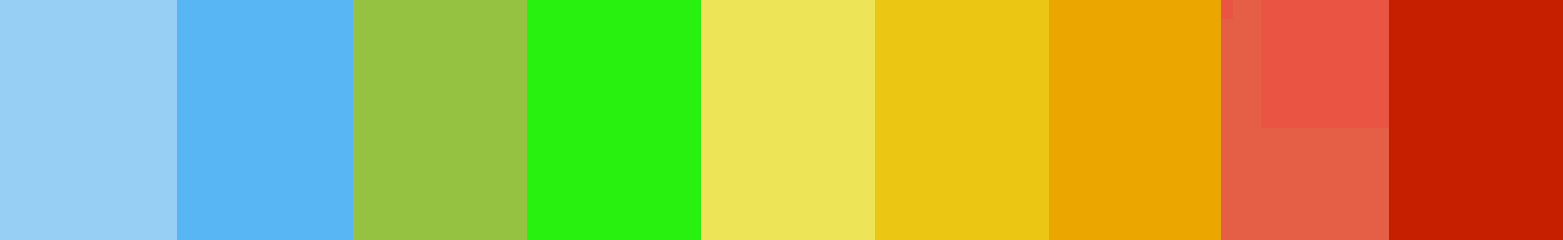}}
\vspace{-0.2cm}
\tiny{\[\hspace{0.7cm}0\hspace{3cm}0.15\text{ cars/m}\]}
\caption{\label{fig:mapMacroscopicModel} A sample output of the density of cars that is calculated using the Payne-Whitham macroscopic model.}
\end{figure}

An example of the output of such a model can be seen in Figure \ref{fig:mapMacroscopicModel}, which shows the traffic density on each segment of the road. 

Our model is derived from the one in \cite{KotsialosPapageorgiou}, the main novelty being that we considered that some of the road segments are controlled by traffic signals. In particular, we focused on the way drivers react to the change of state of a traffic light. We also improved the model in \cite{KotsialosPapageorgiou} by taking into account that there is a speed reduction when cars enter a junction and turn either left or right.

In the following, we present the Payne-Whitham model and outline the particularities that occur at junctions with or without traffic signals. First, we introduce the well-known concept of the ``Fundamental diagram of traffic flow''.

\subsection{Fundamental diagram of traffic flow}
\begin{figure}
\begin{center}\includegraphics[width=\textwidth]{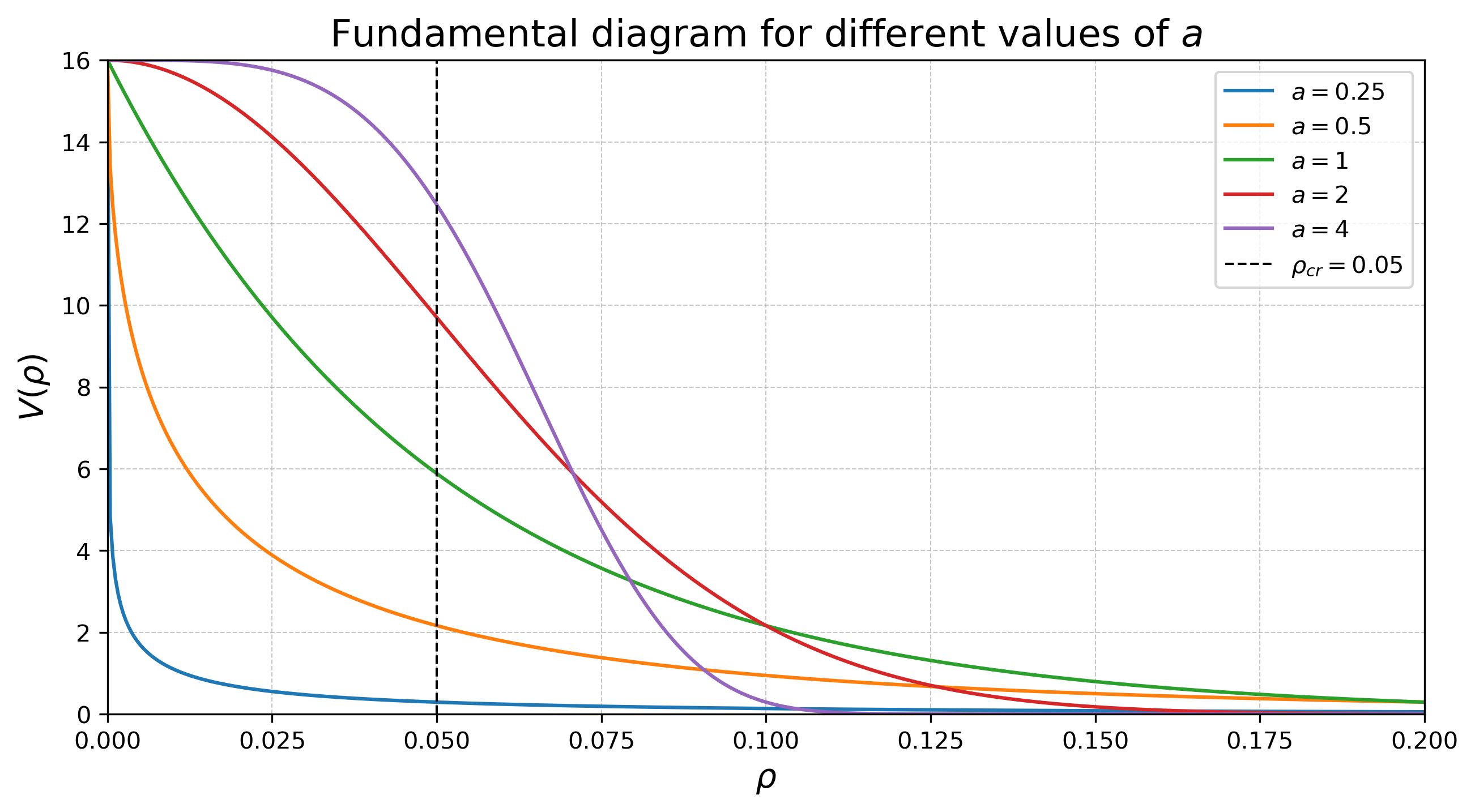}
\caption{\label{fig:fundamentalDiagram} The fundamental diagram of traffic flow for $\rho_{cr}=0.05$ cars/m, \\$v_{max}=16$  m/s  and different values of the parameter $a$.}
\end{center}
\end{figure}

It is an intuitive fact that as the density of cars on a road increases, the speed at which they travel begins to decrease. For any given road, there is a maximum speed $v_{max}$ at which cars usually travel when the road is empty, and also a density $\rho_{cr}$ above which traffic becomes really congested. In addition, the rate at which the speed decays after the critical density is reached may vary depending on the characteristics of the road being analysed. Therefore, when modelling the ideal case of cars driving on a road, in addition to $v_{max}$ and $\rho_{cr}$, we introduce an extra parameter $a>0$ which describes the abruptness of the speed reduction. The ideal relationship between density and speed on a road segment is called the ``fundamental diagram of traffic flow''. In the following, we present the version of the fundamental diagram that appears in \cite[Section II.C]{KotsialosPapageorgiou}:

\begin{equation}
\label{eq:fundamentalDiagram}
V(\rho)=v_{max}\exp\left(-\frac{1}{a}\left(\frac{\rho}{\rho_{cr}}\right)^a\right)
\end{equation}
Here, $\rho$ is the density and $V(\rho)$ represents the theoretical relationship between density and speed on a road segment. This relationship captures the variation in speed with respect to density based on straightforward intuitive principles, such as the observation that vehicles move slowly in congested traffic, but does not explicitly account for interactions between vehicles. Figure \ref{fig:fundamentalDiagram} illustrates the fundamental diagram $V$ for different values of the parameter $a$.


\subsection{A single road segment}
The Payne-Whitham model on a single road segment \cite{Payne1971,Whitham1999} consists of a conservation law for density, together with a characterisation of the instantaneous variation of speed based on the behaviour of drivers in terms of their adaptation to the current density and also their ability to anticipate the changes in density upstream.
\begin{align}
\label{eq:density_continuouos}
&\frac{\partial \rho(t,x)}{\partial t} +\frac{\partial(v(t,x)\rho(t,x))}{\partial x}=0\\
\label{eq:speed_continuouos}
&\frac{\partial v(t,x)}{\partial t}+v(t,x)\frac{\partial(v(t,x))}{\partial x} = \frac{V(\rho(t,x))-v(t,x)}{\nu} - \frac{C}{\rho(t,x)} \frac{\partial \rho(t,x)}{\partial x}.
\end{align}

The meaning of the conservation law \eqref{eq:density_continuouos} is that the density of cars propagates with the speed given by $v(t,x)$.  Similarly, the speed also propagates with the traffic flow, but its variation also depends on other factors given on the right-hand side (RHS) of \eqref{eq:speed_continuouos}. In particular, the first term on the RHS of \eqref{eq:speed_continuouos} drives the speed towards the ideal speed given by the fundamental diagram \eqref{eq:fundamentalDiagram}. The strength of this adaptive behaviour is controlled by the parameter $\nu>0$, which must be calibrated to describe the real-world situation of interest. The second term on the RHS of \eqref{eq:speed_continuouos} represents the ability of drivers to anticipate the traffic conditions ahead and adapt their speed to the change in density. The strength of this ability to anticipate is described by the parameter $C>0$.

\section{Finite difference scheme for the Payne-Whitham model}\label{sec:finite}
The aim of this section is to construct a numerical approximation for the Payne-Whitham model \eqref{eq:density_continuouos} - \eqref{eq:speed_continuouos} on a single road and to provide a stability condition for the time and space discretisation parameters $\delta t$ and $\delta x$. We use an explicit finite difference scheme to compute the density $\rho_e$ and speed $v_e$ at the next time step $k+1$ with respect to the current state (at time step $k$) for each road $e$ uniformly discretised into $N_e$ segments indexed by $i\in \{1,2,\ldots,N_e\}$ (where vehicles travel from the segment $i-1$ directly to the segment $i$):
\begin{equation}
\label{eq:density_numerical}
\frac{\rho_e(k+1,i)-\rho_e(k,i)}{\delta t}+\frac{v_e(k,i)\rho_e(k,i)-v_e(k,i-1)\rho_e(k,i-1)}{\delta x}=0;
\end{equation}
\begin{equation}
\label{eq:speed_numerical}
\begin{aligned}&\frac{v_e(k+1,i)-v_e(k,i)}{\delta t}+\frac{v_e(k,i)^2-v_e(k,i-1)^2}{2\delta x}\\
&\quad\quad=\frac{V(\rho_e(k,i+1))-v_e(k,i)}{\nu}-\frac{C}{\rho_e(k,i)+\chi}\frac{\rho_e(k,i+1)-\rho_e(k,i)}{\delta x}.
\end{aligned}
\end{equation}
Here $\chi>0$ is a parameter used to avoid blow-ups in the last term of \eqref{eq:speed_numerical}. It may also need to be calibrated to reflect the local situation.

\begin{remark}
    On the left hand side of \eqref{eq:speed_numerical} we have chosen to discretise the term: 
    \[v_e(t,x)\frac{\partial(v_e(t,x))}{\partial x} = \frac{\partial((v_e(t,x))^2)}{2 \partial x}\hspace{0.2cm}\text{ numerically as }\hspace{0.2cm}
    \frac{v_e(k,i)^2-v_e(k,i-1)^2}{2\delta x},\]
    due to the increased stability of the model that was practically observed during the experiments. For similar reasons, we have also discretised $V(\rho_e(t,x))$ as $V(\rho_e(k,i+1))$ on the RHS of \eqref{eq:speed_numerical}.
\end{remark}

\begin{remark}
The stability condition to be imposed on the discretisation parameters $\delta t$ and $\delta x$ is of the Courant-Friedrichs-Lewy type (\cite[Section 1.6]{Strikwerda}) and has the form
\begin{equation}
\label{cond:stability}
v_{max}\,\frac{\delta t}{\delta x}\leq 1
\end{equation}

In the case of constant speed, the proof of the stability for the conservation law \eqref{eq:density_numerical} under this condition can be found in \cite[Section 1.6]{Strikwerda}. Furthermore, the suitability of the condition \eqref{cond:stability} for the stability of the Payne-Whitham model was tested in \cite{CaligarisSaconeSiri}.
\end{remark}

\begin{remark}
    We note that in the discrete equations \eqref{eq:density_numerical} and \eqref{eq:speed_numerical}, the upwind spatial finite difference operator was used for the anticipation term (the last fraction in \eqref{eq:speed_numerical}), because anticipating traffic conditions requires information located in front of the driver. 
    
    On the other hand, the downwind finite difference discretisation was applied for the convection terms on the left-hand sides of both equations, since those terms quantify the movement of cars in the direction of the flow.
\end{remark}

\section{Coupling at intersections}\label{sec:intersections}

\subsection{Without traffic lights} 
In the case of a road network, the number of lanes on each street might differ and the cars need to adapt to the change of lane number by either increasing or decreasing the speed. Therefore, the network version of the Payne-Whitham model should take into account the number of lanes $\ell(e)$ corresponding to each edge $e$ of the graph. The next step is to write the equation \eqref{eq:density_numerical} in terms of the traffic flux (i.e. the number of passing cars on all lanes at any time):
\begin{equation}
\label{eq:flux_numerical}
q_e(k,i)\coloneqq\rho_e(k,i)v_e(k,i)\ell(e),
\end{equation}
where by $\rho_e(k,i)$ we mean the density per lane averaged over all lanes of $e$. Therefore, the equation \eqref{eq:density_numerical} can be rewritten as:
\begin{equation}
\label{eq:density_numerical_flux}
\frac{\rho_e(k+1,i)-\rho_e(k,i)}{\delta t}+\frac{1}{\ell(e)}\frac{q_e(k,i)-q_e(k,i-1)}{\delta x}=0.
\end{equation}
Further, it can easily be observed that some of the quantities needed to perform the iteration in \eqref{eq:density_numerical_flux} and \eqref{eq:speed_numerical} are not available from the previous time step. Namely, in order to compute the traffic values (density and speed) for both ends of a road segment $e$ (parameterised by $\{1,2,\ldots, N_e\}$), one would need some virtual values for $q_e(k,0)$, $\rho_e(k,N_e+1)$ and $v_e(k,0)$.

In the case of no traffic lights, we follow the ideas in \cite[Section II.D]{KotsialosPapageorgiou}, where these virtual values are defined as weighted sums of the traffic values on the other road segments adjacent to this intersection. The following three formulas were adapted from \cite{KotsialosPapageorgiou} in order to prepare the path for integrating traffic signals. We also attempted to slightly improve them according to our intuition about modeling traffic in intersections and, moreover, to correct several misprints.

\begin{figure}
    \centering
    \includegraphics[scale=0.36]{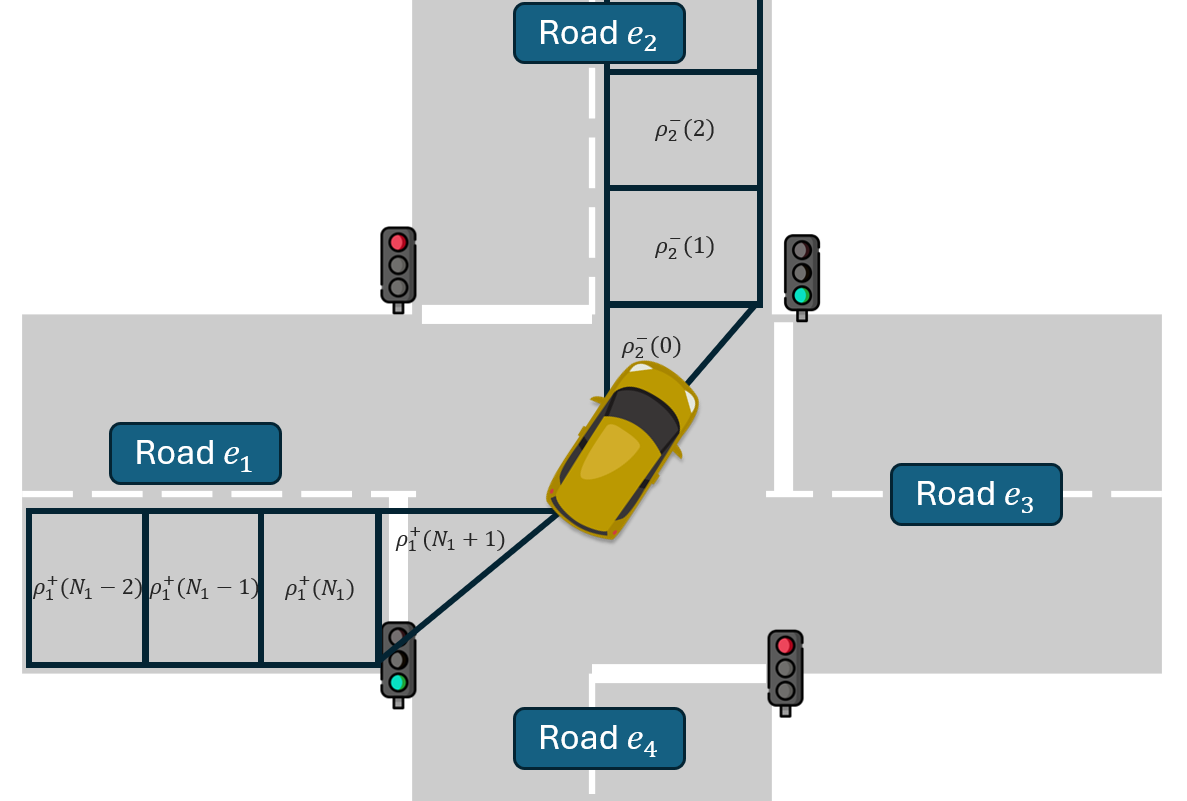}
    \caption{An intersection with the virtual initial and final densities, $\rho_2^{-}(0)$ and $\rho_{1}^{+}(N_1+1)$, respectively. The signs $\pm$ in the notation are used to help differentiate between the incoming and outgoing sides of each street, in relation to the depicted intersection.}
    \label{fig:intersection1}
\end{figure}

We illustrate these calculations for the particular case of the intersection in Figure \ref{fig:intersection1} (see the caption for notations):
\begin{align}
\label{eq:density_virtual_initial}
{q_{2}^-(k,0)} &=\sum_{i=1}^4 q_{i}^+(k,N_{i})\cdot \frac{\omega(e_i,e_2)}{\sum_{j=1}^4 \omega(e_i,e_j)};\\
\label{eq:speed_virtual_initial}
{v_{2}^-(k,0)} &=\frac{1}{q_{2}^-(k,0)}\sum_{i=1}^4 v_{i}^+(k,N_{i})\cdot q_{i}^+(k,N_{i})\cdot \frac{\omega(e_i,e_2)}{\sum_{j=1}^4 \omega(e_i,e_j)};\\
\label{eq:density_virtual_final}
{\rho_{1}^+(k,N_{1}+1)} & =\frac{\sum_{i=1}^4 (\rho_{i}^-(k,1))^2\cdot \omega(e_1,e_i)}{\sum_{i=1}^4 \rho_{i}^-(k,1)\cdot \omega(e_1,e_i)}.
\end{align}

The turning weights $\omega(e_i,e_j)$ represent the number of cars that, over a large period of time, choose to continue their trip on road $e_j$, after entering that particular intersection from road $e_i$. These weights encode the average preference of drivers for following a certain path in the intersection and are determined by processing real-world traffic data.
 
The intuition behind \eqref{eq:density_virtual_initial} is that the virtual flux at the entrance of a street from a junction is the sum of the fluxes on each road entering the intersection, weighted with the rate of the vehicles on that incoming road that will actually enter that particular street.

Formula \eqref{eq:speed_virtual_initial} has a similar interpretation. Here, the incoming speeds are weighted based on the actual number of cars (i.e. flux) that has that speed. The total flux through which we divide is the one that reaches the particular street, so it is the same as in the above formula \eqref{eq:density_virtual_initial}, whereas the relative flux $ q_{i}^+(k,N_{i})\cdot \frac{\omega(e_i,e_2)}{\sum_{j=1}^4 \omega(e_i,e_j)}$ represents the number of cars that will travel from $e_i$ to $e_2$ in a unit of time.
 
Eventually, formula \eqref{eq:density_virtual_final} is equivalent to:
\[{\rho_{1}^+(k,N_{1}+1)}  =\frac{\sum_{i=1}^4 \rho_{i}^-(k,1) \cdot \rho_{i}^-(k,1) \cdot \frac{\omega(e_1,e_i)}{\sum_{j=1}^4\omega(e_1,e_j)}}{\sum_{i=1}^4 \rho_{i}^-(k,1)\cdot \frac{\omega(e_1,e_i)}{{\sum_{j=1}^4\omega(e_1,e_j)}}},\]
which has the following meaning: the virtual density felt by drivers when entering an intersection from one street is influenced by the densities they will encounter when entering the outgoing roads, weighted by the density of cars that choose to enter each of those outgoing roads. The heavier the traffic on one outgoing street, the more significant that density will obstruct the drivers from entering the intersection towards that street.

Furthermore, we think that the explanation for the quadratic term at the numerator of \eqref{eq:density_virtual_final} given in \cite{KotsialosPapageorgiou} is valid:
``The quadratic
term is used in \eqref{eq:density_virtual_final} to account for the fact that one congested
leaving link may block the entering link even if there is free flow
in the other leaving link.''


\subsection{Limiting the turning speed}
Our first improvement to the algorithm in \cite{KotsialosPapageorgiou} concerns the maximum speed of cars changing direction at an intersection. More precisely, we have limited the virtual initial speed $v_e(k,0)$ according to the speed at the other edges in the intersection and the geometry of the intersection. As usual, we illustrate this procedure using the notation in Figure \ref{fig:intersection1}: the speed value given by \eqref{eq:speed_virtual_initial} must not exceed the following limit:
\begin{align}
\label{eq:speed_limit_initial}
&v^{lim}_{2}(k,0)\coloneqq\frac{1}{q_{2}^-(k,0)}\sum_{i=1}^4 v^{max}_{i}\frac{(1-\cos(\widehat{e_i,e_2}))}{2}q_i^+(k,N_{i})\cdot \frac{\omega(e_i,e_2)}{\sum_{j=1}^4 \omega(e_i,e_j)},
\end{align}
where we recall that $v^{max}_{i}$ is the speed at which cars normally travel when the road $e_i$ is empty and by $\widehat{e_i,e_j}$ we understand the angle between the edges $e_i$ and $e_j$ in the traffic network. 

For instance, in the case of going straight (i.e. $\widehat{e_i,e_j}=180^\circ$) the cars are allowed to travel at maximum speed, but if they make a $90^\circ$ turn, the speed is reduced to half of $v_{max}$. Finally, in the case of a U-turn (i.e. an angle of $360^\circ$), the cars must stop before the turn.

\subsection{Adding traffic lights}

Our major contribution to the traffic simulation model is that we consider some of the intersections being controlled by traffic signals. The modifications that occur in this situation imply both the speed values at the end of the roads that are controlled by traffic lights and the weighted sums \eqref{eq:density_virtual_initial}-\eqref{eq:speed_limit_initial}.

If the traffic signal is {\color{red} RED}:
\begin{itemize}
\item The final speed $v_e(N_e)$ of the edge $e$ is set to zero;
\item The road $e$ is not taken into account as an input in the weighted sums for the virtual density and speed of any other edge. More precisely, we set $\omega(e,e_i)=0$, for each road $e_i$.
\end{itemize}


\section{Initialising and calibrating the model with real-world traffic data}\label{sec:init}
In this section, we describe how we prepared the experimental scenario, including the information about the road network and real-world traffic data that we used, and how we have integrated the data to both calibrate and initialise our model. 


\subsection{Road network and intersection data}

The aim of this section is to explain the data used to describe the road structure and the behaviour of drivers at intersections. First, the road network was extracted from the open geographic database, OpenStreetMap (OSM) \cite{openstreetmap2}. In particular, we limited our research to a central district of the city of Warsaw, in Poland, namely Stara Ochota, with $21$ intersections with traffic lights, and we requested information about the road segments, their coordinates, the type of roads, the location of the intersections and the list of their incoming and outgoing roads, and the locations of the traffic lights, also indicating for which road entry they are responsible. The choice of such a district (which is a sub-region of the whole road network of Warsaw) was motivated by the fact that in some previous experiments, it turned out that considering the whole road network of Warsaw can be computationally demanding \cite{gora2011genetic}. On the other hand, the considered area should contain at least several intersections to ensure that the coordination of traffic signals between different intersections and the application of our optimisation algorithms can bring some benefits. Also, the same area was used in several other research works (e.g., \cite{skowronek2021graph, skowronek2023graph, gora2019solving}). Finally, since the OSM service data sometimes has limited correspondence with the real road network, the central area of Warsaw was chosen because it has relatively good quality OSM data compared to, for example, the suburban areas.

We also used the data describing the probabilities of drivers following each particular path at every intersection. It was obtained using the Traffic Simulation Framework \cite{gora2009traffic}. Specifically, the extracted information consists of the turning fraction $\omega(e_i,e_j)$ for each road segment $e_i,e_j$ in every intersection. 

The diagram \ref{diagram:osm} summarizes this process.

\begin{tcolorbox}[colframe=black,colback=white,boxrule=1pt]
\vspace{-1.2em}
\begin{equation}\label{diagram:osm}
\begin{aligned}
\text{OpenStreetMap Data }&\Rightarrow
\left\{\begin{aligned} &\text{ The road structure: }(e_i)_i,\, \ell(e_i),\, \widehat{e_i,e_j};\\
&\text{ The location of traffic signals};
\end{aligned}\right.\\
\text{Traffic Simulation Framework } &\Rightarrow \, \omega(e_i,e_j).
\end{aligned}
\end{equation}
\end{tcolorbox}

\subsection{Calibration of the fundamental diagram}\label{sec:calibrate-fundamental}

The first step in the calibration process requires choosing the correct parameters of the fundamental diagram of traffic \eqref{eq:fundamentalDiagram}. In other words, we need to approximate as accurately as possible the real dependence between density and speed in urban traffic.
Luckily, with the help of the Traffic Technical Institute at the Faculty of Civil and Geodetic Engineering, University of Ljubljana, we obtained data from traffic counters at several locations in Ljubljana, Slovenia. More precisely, the datasets that we received contain the number of cars together with their average speeds, values aggregated at 5-minute intervals. The density of cars was computed from the available data using the following formula:
\[{\rm density}(cars/m) = \frac{{\rm number\_of\_cars}(cars)}{{\rm timeframe}(s)\cdot {\rm average\_speed}(m/s)}.\]

Eventually, with the empirical density and speed measurements at hand, we have fitted the three parameters of the fundamental diagram \eqref{eq:fundamentalDiagram} using the functionality \code{curve\_fit} from the Python library \code{scipy}\cite{2020SciPy-NMeth}, thus obtaining the following parameters:
\[v_{max}= 13.68\, m/s, \quad \rho_{cr}= 0.05\, cars/m, \quad a= 1.24.\]
 The latter two of those parameters will be used for every road segment in the model, whereas the former, $v_{max}$, will be replaced by the free flow speed retrieved during the initialisation process described in the next section.

The diagram \ref{diagram:calibration} summarizes the calibration process.

\begin{tcolorbox}[colframe=black,colback=white,boxrule=1pt]
\vspace{-1.2em}
\begin{align}\label{diagram:calibration}
&\text{Traffic data from car counters in Ljubljana } \nonumber \\&\Rightarrow \text{Calibration of }v_{max},\, \rho_{cr},\, a
  \text{ in the fundamental diagram (see \eqref{eq:fundamentalDiagram})}.
\end{align}
\end{tcolorbox}


\subsection{Initialisation of the model} 
\label{sec:initial-data}
The aim of this step is to provide the initial data $(\rho_e(0,i),v_e(0,i))$, $i\in\{1,2,\ldots,N_e\}$ on each edge of the network, corresponding to the time step $k=0$ of the evolution equations \eqref{eq:density_numerical}-\eqref{eq:speed_numerical}. The initial data should reflect the real traffic conditions in the city we analyse.

In order to do this, we have acquired traffic data (average speeds, maximum (free flow) speeds, and travel times) from TomTom's Flow Segment Data API \cite{tomtom} concerning the Stara Ochota district in Warsaw.

According to this data, we have overwritten the free flow speed of the fundamental diagram (i.e. $v_{max}$ in Section \ref{sec:calibrate-fundamental}) with the free flow speed provided by TomTom's API. The other parameters ($\rho_{cr}$ and $a$) were left unchanged from Section \ref{sec:calibrate-fundamental}, since the data from TomTom contained only information regarding vehicles' speed.

Next, since we need to initialise each of the $N_e$ segments in the discretisation of each road, the spatial resolution of the real-world data used in the initialisation process should be very fine. Therefore, we divided the large road segments in TomTom's data representation into smaller parts (of approximately 10 meters each) and introduced small random deviations from the unique value of the average speed on the large segment, returned by the API. Thus, we obtained the speed $v_e(0,i)$ for each 10-meter road division $i\in\{1,2,\ldots,N_e\}$ and then we inverted the fundamental diagram to obtain the corresponding density: \[\rho_e(0,i)\coloneqq V^{-1}(v_e(0,i)),\, \forall i\in\{1,2,\ldots,N_e\}.\]
Eventually, we used the obtained pairs $(\rho_e(0,i),v_e(0,i))$ for each road division to initialise the model.

The diagram \ref{diagram:init} summarizes the initialisation process.

\begin{tcolorbox}[colframe=black,colback=white,boxrule=1pt]
\vspace{-1.2em}
\begin{equation}\label{diagram:init}
\text{TomTom's API }\Rightarrow \left\{
\begin{aligned}
&\text{override }v_{max}\text{ from \eqref{diagram:calibration}};\\
&\text{the intial data of \eqref{eq:density_numerical}-\eqref{eq:speed_numerical}: }(\rho_e(0,i),v_e(0,i)).
\end{aligned}
\right.
\end{equation}
\end{tcolorbox}

\subsection{Setting the values of the global parameters}
Our model uses several global parameters that should be chosen to reflect realistic traffic behaviour, especially near traffic signals. More specifically, these parameters are $\nu$, $C$ and $\chi$ from \eqref{eq:speed_numerical} and they are fitted in order to obtain a realistic behaviour of cars in the presence of traffic lights as we observe the model run multiple times. After considering several combinations of parameters, we have chosen them such that the behaviours they describe (i.e., the capability of drivers to adapt to the fundamental diagram and to anticipate the traffic conditions ahead) match a realistic formation of a queue in front of a \textcolor{red}{RED} signal. Furthermore, the chosen parameters also facilitate the dispersion of the high-congestion area behind the signal after the traffic light turns \textcolor{Green}{GREEN}.

\begin{figure}[h]
\begin{center}
\vspace{-0.7cm}
\includegraphics[scale=0.5]{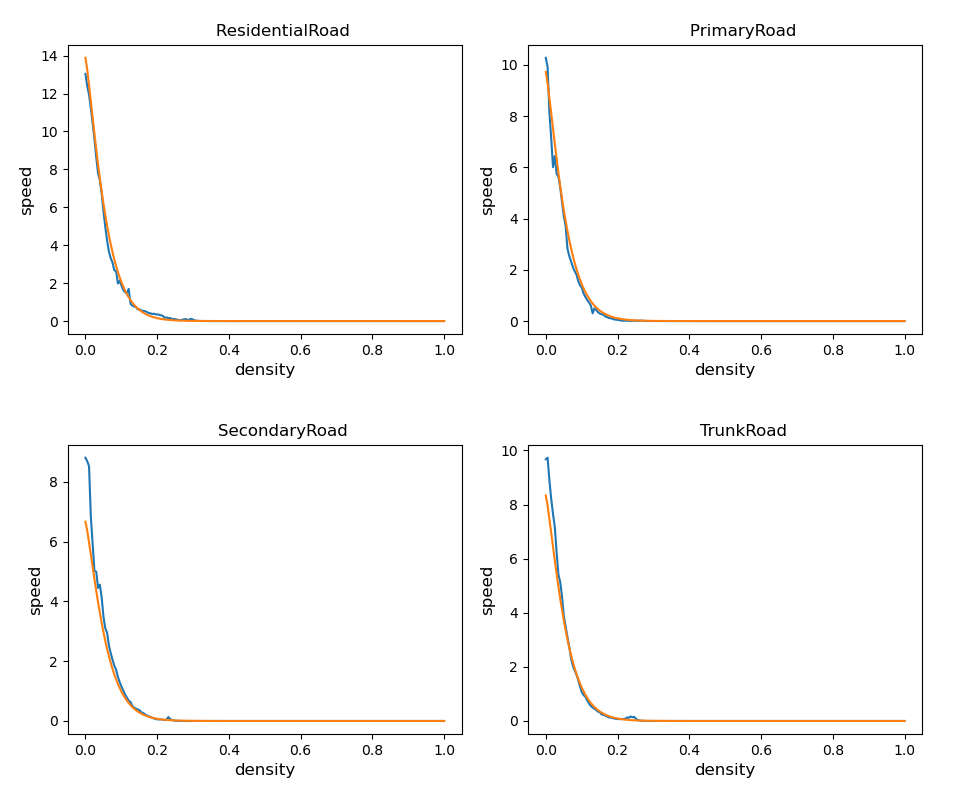}
\caption{\label{fig:numerical_theoretical_fundamental_diagram} The relation between density and speed obtained by running the calibrated model (blue) is consistent with the fundamental diagram (orange) for all four categories of roads.
}
\end{center}
\end{figure}

In addition to that, we have fine-tuned the parameters such that the relationship between density and speed is realistic and consistent with the fundamental diagram \eqref{eq:fundamentalDiagram}. Figure \ref{fig:numerical_theoretical_fundamental_diagram} shows the comparison between simulated values for density and speed and the values given by the fundamental diagram \eqref{eq:fundamentalDiagram} for four categories of roads, after 180 seconds of simulated traffic. The speed is slightly lower or higher than in the ideal case, probably because of vehicles that are stopped at a \textcolor{red}{RED} signal or turn from one street to another. 

Eventually, the fitted values of the parameters of the model are:
\[\nu=1, \quad C=\, 7\, m^2/s, \quad \chi=0.008\, cars/m.\]

The diagram \ref{eq:setting} summarizes the choice of these parameters.

\begin{tcolorbox}[colframe=black,colback=white,boxrule=1pt]
\vspace{-1.2em}
\begin{equation}\label{eq:setting}
\begin{aligned}
    \text{Observations of multiple }\hspace{1em}& \\ \text{runs of the PWTL model }\hspace{1em}&
    \end{aligned}\Rightarrow \begin{aligned}\hspace{1em}\nu,\, C,\, \chi &\text{ for a realistic behaviour}\\
    &\text{ near traffic lights.}
    \end{aligned}
\end{equation}
\end{tcolorbox}

\section{Optimization procedure}\label{sec:optim}

In extending the Payne-Whitham model to traffic lights, our aim was to create a framework that would allow the simulation of multiple traffic light configurations in order to find the optimal one. Sub-optimal traffic light configurations are often responsible for traffic congestion, which leads to increased time spent in traffic and consequently increased vehicle emissions.

As presented in Section \ref{sec:intro}, there have been many approaches to optimising the configuration of traffic lights. Some of them are based on traffic simulation frameworks (e.g. \cite{gora2019solving, wei2019colight}), which are used to simulate traffic conditions in different configurations and select the optimal one. This selection is made using an optimisation technique, such as genetic algorithms \cite{gora2019solving}, or machine learning methods, such as reinforcement learning \cite{wei2019colight}.

A notable challenge in optimising traffic light configurations through simulation is its inherent computational and time-intensive nature \cite{ismis}. Given the large space of possible configurations, this computational complexity becomes a significant concern. To address this, we turn our attention to surrogate models.

\subsection{Surrogate modeling}

Surrogate models serve as efficient approximations of computationally expensive traffic simulation models \cite{forrester2008engineering}. By capturing the essential relationships between input parameters and simulation results, surrogate models significantly reduce computational costs. This allows for faster exploration of the vast configuration space, making the optimisation process more tractable and feasible for large-scale simulations.

In this research, we have chosen to use a set of different possible surrogate models, with different levels of complexity, as described in Section \ref{sec:sme}. These surrogate models are designed to approximate the output of the PWTL model in a computationally cheaper way, by learning the relationship between traffic light configurations and key congestion variables derived from multiple simulated traffic conditions. To achieve this, we first conduct a series of simulations using the PWTL model, each with specific initial traffic light configurations. From these simulations, we extract relevant congestion metrics. The surrogate models are then trained in a supervised manner using k-fold cross-validation to generalize this relationship, enabling efficient predictions without the need for additional simulations.
 

\subsection{Traffic light configuration space}

In our approach, the traffic light configuration space is characterised by three key variables associated with each intersection in our dataset: the red light duration, the green light duration, and the time elapsed from the start of the simulation to the first red-to-green light change, henceforth referred to as the ``offset time''. They encapsulate the temporal aspects of traffic signal control that affect the flow of vehicles through intersections. These variables are set at the beginning of the simulation and do not change during the whole simulation.  

\begin{figure}
\begin{center}
\includegraphics[scale=0.46]{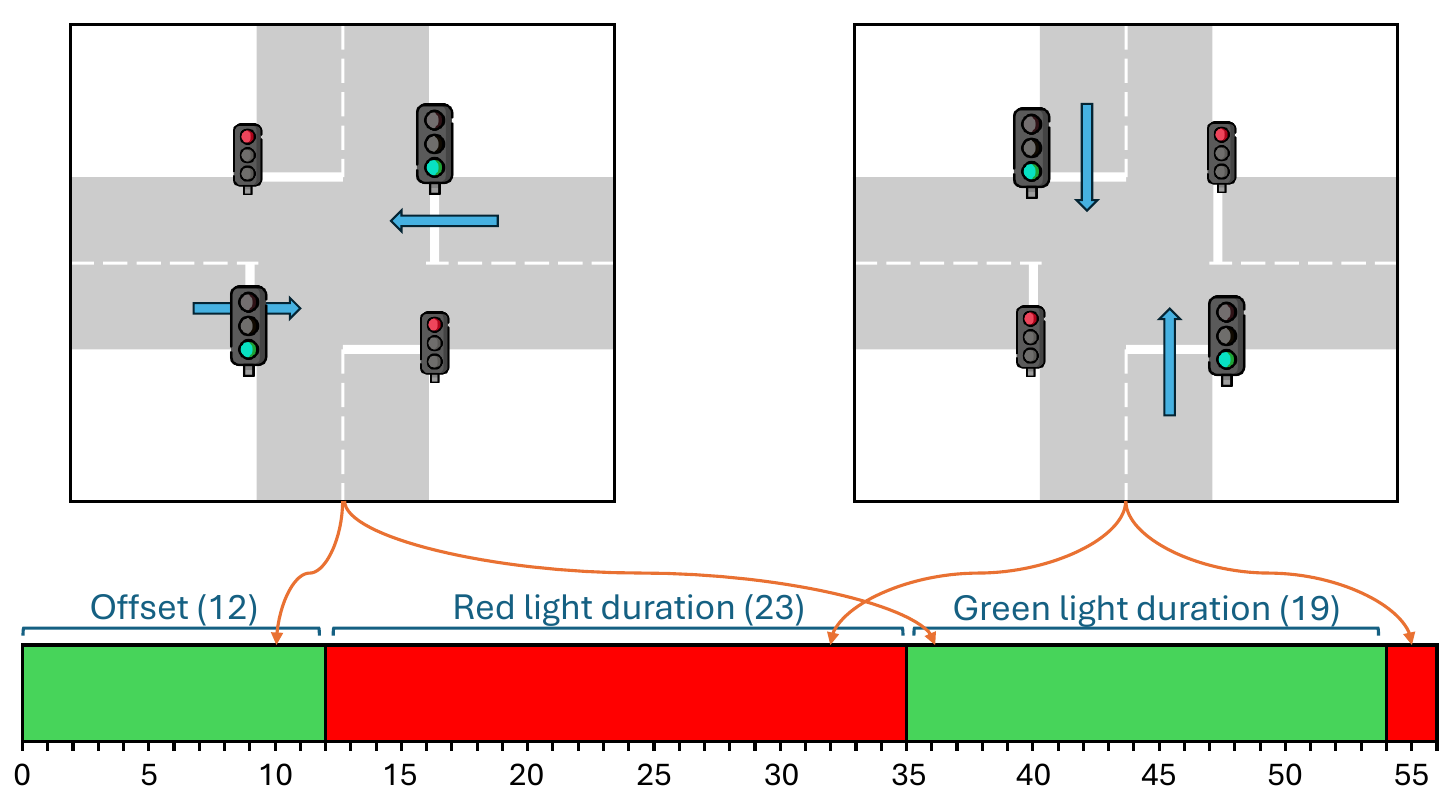}
\caption{\label{fig:intersection2} An example of a 4-way intersection, where opposing road segments share the same traffic light configuration. The diagram highlights three key traffic light variables: an offset of 12, a red light duration of 23, and a green light duration of 19. In this example, the signal is \textcolor{Green}{GREEN} on the horizontal direction at the beginning of the simulation.}
\end{center}
\end{figure}

In our modelling approach, we made a specific assumption about intersections equipped with traffic lights. We assumed that for each such intersection there are exactly two sets of synchronised traffic lights. For example, at an intersection with four roads, each set of opposing roads shares an identical light configuration, as shown in Figure \ref{fig:intersection2}. This synchronisation simplifies the modelling process and is in line with common traffic engineering practice, where opposite directions often have simultaneous green lights to facilitate smoother traffic flow and reduce potential conflicts.

\subsection{Objective function}

\begin{figure}
\begin{center}
\includegraphics[scale=0.7]{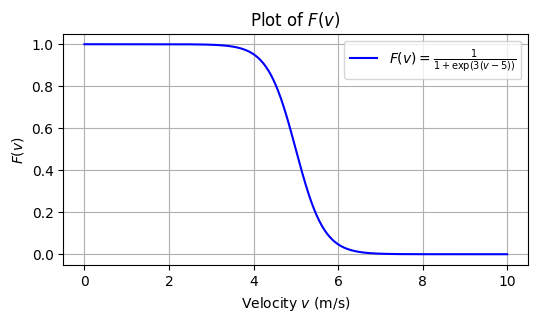}
\caption{\label{fig:queue_f} The representation of the logistic function used to assess whether a road segment is congested or not. In particular, a segment is considered congested if the speed falls below $5\ m/s$.}
\end{center}
\end{figure}

When running a simulation using the PWTL framework, we focus on assessing the level of congestion within the road network. This assessment is based on two key variables measured at each time step of the simulation:

\begin{itemize}
\item \textbf{Average Speed}: This variable represents the average speed of vehicles across all roads in the network. It serves as an indicator of the overall traffic flow efficiency.
\item \textbf{Queue Length}: To measure the length of queues formed along the roads, we first define when a road segment is considered congested. A segment is classified as congested if its speed falls below a specified threshold (hereafter referred to as $v_q$). To smooth the transition between congested and non-congested states, we apply a logistic function $F$ to the speed $v$ of the road segment. The function $F$ is defined as follows:
 \[F(v)=\frac{1}{1+\exp(c(v-v_{q}))}\]   
In our simulation, we set $c=3$, and $v_{q}=5\ m/s$. Figure \ref{fig:queue_f} illustrates the logistic function $F$ using these parameters. This function allows for a gradual classification of congestion, rather than a binary threshold. To measure the total queue length of the network, we sum the road congestion $F(v)$ for all divisions (segments of approximately $10$ meters) of the entire road network.
\end{itemize}

For each simulation, we take the average of these variables over all time steps to provide a comprehensive measure of traffic congestion throughout the simulation period.

\subsection{Optimization technique: Differential Evolution}

After having trained the surrogate model, we use an optimization technique to find the optimal traffic light configuration that either maximizes the average speed or minimizes the queue length. To efficiently navigate the large parameter space, we employ the Differential Evolution (DE) algorithm \cite{storn1997differential}. This optimization algorithm is particularly well suited to our traffic light optimization task due to its robustness and effectiveness in handling complex, multi-modal, and high-dimensional optimization problems. The DE algorithm excels in exploring the global search space and avoiding local minima, making it an ideal choice for our application. The search for the optimum configuration begins with the best traffic light setup identified during the training phase.

\section{Optimization experiments}\label{sec:optexp}

In the following sections, the setup and results of our optimization experiments are illustrated.

\subsection{Traffic Simulation}

In order to generate a comprehensive dataset for the training of our surrogate models, we adopted the following approach. The input data generation strategy involved simulating traffic using the PWTL model described above for a fixed duration of 340 in-game seconds using the same initial traffic conditions, as described in Section \ref{sec:initial-data}. To ensure the results depend less on the initial traffic conditions and allow the model to stabilize, we started measuring the traffic congestion variables only after the first 100 in-game seconds. During each simulation run, a traffic light configuration was randomly selected from the predefined set of possible combinations for each intersection. This set is the one that contains all combinations for the following configurations:
\begin{itemize}
            \item Red light duration (in in-game seconds) $\in \{20,21,22,\dots,53,54\}$
            \item Green light duration (in in-game seconds) $\in \{20,21,22,\dots,53,54\}$
            \item Offset time (in in-game seconds) $\in \{0,1,2,\dots,\text{Red light}+\text{Green light} -1 \}$  
        \end{itemize}
The average speed of the vehicles and the number of congested road segments in the whole simulated road network over the whole period was then calculated as the output. 

\begin{figure}
\begin{center}
\includegraphics[width=\textwidth]{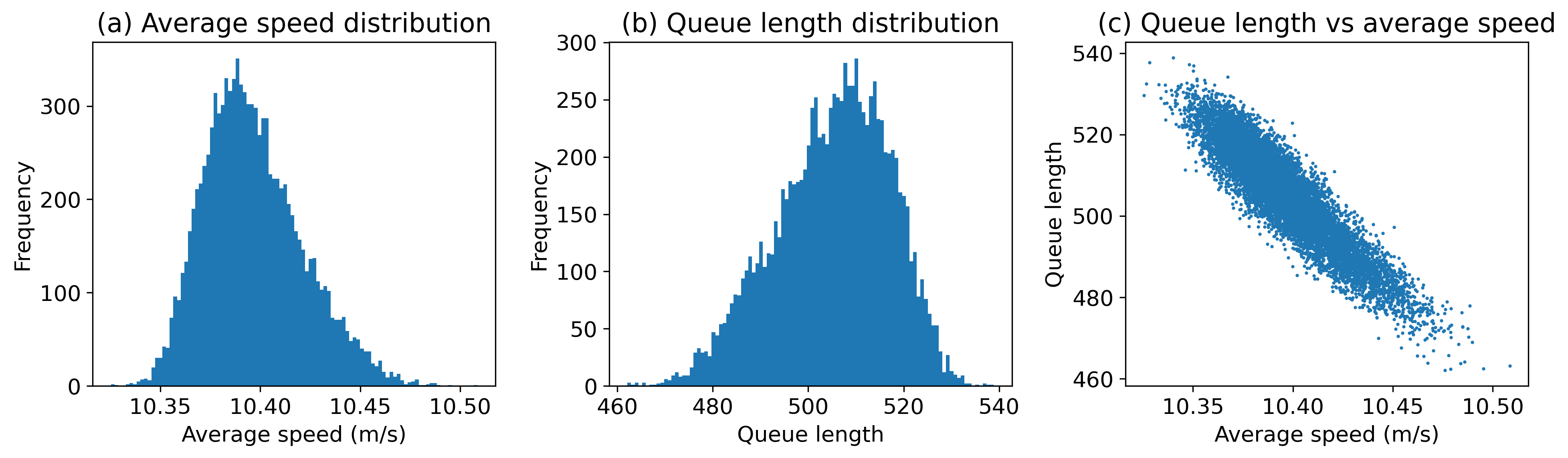}
\caption{\label{fig:data_distr} Distribution of the traffic congestion variables outputted by the PWTL framework. Panel (a) shows the distribution of the average speed, which follows a bell-curve distribution, slightly skewed on the right side. Panel (b) shows the distribution of the queue length variable, which also follows a skewed bell-curve distribution. Panel (c) shows the scatter plot of all simulations, with the average speed on the $x$-axis and the corresponding queue length on the $y$-axis. It's apparent that there is a negative relationship between the two variables, confirming the intuitive notion that a lower number of congested road segments increases the overall average speed.}
\end{center}
\end{figure}

This process was repeated for 10000 simulations (running a single simulation took about $7$ minutes on a standard CPU \footnote{Intel(R) Xeon(R) CPU E5-2698 v4 @ 2.20GHz}, but to speed up the calculations we used $40$ cores of the computational server ``Chuck'' available at the University of Warsaw), generating a diverse dataset that captured the performance of different traffic light configurations under the same initial conditions. This dataset then provided the training and test inputs for the surrogate models, allowing them to learn the intricate relationships between the chosen traffic light configurations and the resulting traffic congestion variables. Figure \ref{fig:data_distr} shows the distribution of the congestion variables obtained from the traffic simulations.

\subsection{Surrogate models and evaluation}\label{sec:sme}

In the context of this article, we trained three types of surrogate models: one that predicts the average speed, one that predicts the total queue length, and one that predicts both variables simultaneously.

For these three models, we employed Linear Regression, the Multi-Layer Perceptron (MLP) \cite{goodfellow2016deep}, the Random Forest (RF) \cite{breiman2001random}, and the Support Vector Regression (SVR) \cite{drucker1996support} algorithms. Linear Regression provides a straightforward approach to modelling linear relationships between inputs and outputs. MLP, a type of neural network, is well-suited for capturing non-linear relationships due to its multiple layers and activation functions. Random Forest, an ensemble learning method, offers robustness and high accuracy by aggregating the predictions of multiple decision trees. Additionally, we used Support Vector Regression (SVR) for models with one-dimensional outputs, which is effective for regression tasks with smaller datasets and high-dimensional spaces. SVR is not used for multi-output models because it is designed for single-output regression. Creating multi-output models with SVR would require developing two independent models, one for speed and one for queue length, which would be equivalent to the single-output models previously generated. In Appendix \ref{app:surrogate_par}, we describe the set of parameters that were explored for each algorithm.
Additionally, the input data regarding the traffic light configuration for each intersection has been normalized between 0 and 1.

To measure accuracy, we chose to use the Root Mean Squared Error (RMSE). Considering \(x_i\) as one of the output variables, such as the simulated total average speed (from the whole simulation on the whole road network) and \(f_i\) the predicted average speed using the surrogate model, for each simulation \(i \in \{1,2,\dots,n\}\), the RMSE metric is:

\[
\text{RMSE} = \sqrt{\frac{1}{n}\sum_{i=1}^n (x_i - f_i)^2}
\]

The RMSE has the same scale as the original data, making it more interpretable than the Mean Squared Error (MSE). Since RMSE squares errors before averaging and then takes the square root, it places greater weight on larger discrepancies, making it particularly useful when large errors need to be penalized more in model evaluation. Tables \ref{tab:results_speed}, \ref{tab:results_queue}, and \ref{tab:results_both} show the results according to RMSE for all the models considering the best-performing set of parameters. For the models that predict the average speed and the queue length simultaneously, we calculate the average of the RMSE for each variable.

\subsection{Optimization and validation}

In this section, we analyze the results obtained when applying the Differential Evolution optimization algorithm (using the \code{scipy.optimize} library) to the surrogate models to find the traffic light configuration that optimizes the traffic viability. 

For surrogate models predicting either speed or queue length, we optimized the respective variables. For models predicting both, we optimized each of its output
variables independently
, determining separate optimal configurations for speed and for queue length. We then validated these optimized configurations by simulating traffic using our PWTL framework to assess their performance according to the calibrated PWTL model, rather than relying solely on the surrogate models.

\begin{table}[h]
    \centering
    \begin{tabular}{|l|l|l|l|l|}
        \hline
        & \textbf{LR} & \textbf{SVR} & \textbf{MLP} & \textbf{RF} \\ \hline
        \textbf{RMSE} & 0.089 & 0.0112 & 0.0143 & 0.0122 \\ \hline
        \textbf{Predicted max speed} & 10.404 & 10.673 & 10.548 & 10.494 \\ \hline
        \textbf{Simulated avg. speed} & 10.564 & 10.497 & 10.561 & 10.509 \\ \hline
        \textbf{Simulated queue length} & 452.0 & 467.9 & 441.9 & 463.2 \\ \hline
    \end{tabular}
    \vspace{1em}
    \caption{Results related to surrogate models that predict speed. Each column represents a different surrogate model, and the rows showcase the results for the accuracy of the model's prediction (RMSE), the optimal speed (according to the given surrogate model)  found using the Differential Evolution (DE) algorithm, and the results from the PWTL simulation using the optimal traffic light configuration.}
    \label{tab:results_speed}
\end{table}

In Table \ref{tab:results_speed}, we present the outcomes for the surrogate models trained to predict the average speed across the road network. Based on RMSE, SVR emerges as the most accurate surrogate model in replicating the PWTL framework for our dataset, followed by Random Forest and MLP. Linear Regression exhibits notably lower accuracy. Interestingly, SVR predicts the highest speeds among the models, whereas Linear Regression predicts lower speeds. However, during running traffic simulations using optimal traffic light configurations, Linear Regression achieves the highest speeds overall, while MLP results in the shortest queue lengths. Surprisingly, despite its accuracy in training, SVR performs the worst in terms of simulated speed and queue length for the optimal traffic signal settings. For context, within our training data, the maximum speed recorded was $10.509\ m/s$, and the minimum queue length was $462.1$.

Upon investigating the underlying reasons, we conducted a detailed analysis of the specific traffic light configurations. It was observed that Linear Regression consistently favours a configuration strategy where, for each intersection, the red light duration for one direction is set to the maximum allowed time of 54 seconds, while the green light duration is set to the minimum allowed time of 20 seconds. This approach suggests that the LR model optimizes speed by prioritizing traffic flow in one direction over the other at all intersections. In contrast, the traffic light configurations identified by SVR exhibit less distinct patterns.

\begin{table}[h]
    \centering
    \begin{tabular}{|l|l|l|l|l|}
        \hline
        & \textbf{LR} & \textbf{SVR} & \textbf{MLP} & \textbf{RF} \\ \hline
        \textbf{RMSE} & 42.06 & 6.505 & 5.646 & 5.555 \\ \hline
        \textbf{Predicted min queue length} & 501.8 & 411.3 & 379.2 & 466.7 \\ \hline
        \textbf{Simulated avg. speed} & 10.523 & 10.484 & 10.513 & 10.476 \\ \hline
        \textbf{Simulated queue length} & 444.5 & 465.8 & 455.0 & 462.1 \\ \hline
    \end{tabular}
    \vspace{1em}
    \caption{Results related to surrogate models that predict queue length. Each column represents a different surrogate model, and the rows showcase the results for the accuracy of the model's prediction (RMSE), the optimal queue length (according to the given surrogate model) found using the Differential Evolution algorithm, and the results from the PWTL simulation using the optimal traffic light configuration.}
    \label{tab:results_queue}
\end{table}

Table \ref{tab:results_queue} displays the outcomes obtained from surrogate models trained to predict queue length. Once again, the straightforward LR algorithm exhibits the lowest accuracy based on RMSE. However, LR manages to identify the optimal traffic light configuration during simulations, achieving superior performance in both speed and queue length. In contrast, while the RF algorithm demonstrates higher accuracy in predictions, its performance in simulations proves less effective. 

\begin{table}[h]
    \centering
    \begin{tabular}{|l|l|l|l|l|}
        \hline
        & & \textbf{LR} & \textbf{MLP} & \textbf{RF} \\ \hline
        & \textbf{RMSE} & 21.08 & 3.114 & 2.791 \\ \hline 
        \multirow{4}*{\textbf{Opt. for speed}} 
        & \textbf{Predicted max speed} & 10.404 & 10.604 & 10.488 \\ \cline{2-5}
        & \textbf{Corresponding queue length} & 411.2 & 474.7 & 485.5 \\ \cline{2-5}
        & \textbf{Simulated avg. speed} & 10.561 & 10.495 & 10.509 \\
        \cline{2-5}
        & \textbf{Simulated queue length} & 451.5 &	463.8 & 463.2 \\ \hline
        \multirow{4}{*}{\textbf{Opt. for queue length}} 
        & \textbf{Predicted min queue length} & 501.8 & 227.4 & 469.6 \\ \cline{2-5}
        & \textbf{Corresponding speed} & 10.610 & 10.504 & 10.436 \\ \cline{2-5}
        & \textbf{Simulated avg. speed} & 10.523 & 10.522 & 10.509 \\
        \cline{2-5}
        & \textbf{Simulated queue length} & 444.5& 443.5 & 463.2 \\ \hline
    \end{tabular}
    \vspace{1em}
    \caption{Results related to surrogate models that predict average speed and queue length together. The three columns on the right represent different surrogate models, and the rows showcase the results for the accuracy of the model's prediction (RMSE), and the model's optimum together with the corresponding simulation results when using the Differential Evolution algorithm to find the optimal speed and the optimal queue length.}
    \label{tab:results_both}
\end{table}

Table \ref{tab:results_both} showcases the results for surrogate models predicting both average speed and queue length together. Once again, Linear Regression (multiple LR in this case) is the worst-performing model for RMSE, but shows remarkable simulation results both when optimizing for speed and for queue length. MLP and RF, despite being very accurate surrogate models, are not able to predict configurations that are optimal during the traffic simulation.

\begin{figure}[h]
\begin{center}
\includegraphics[scale=0.5]{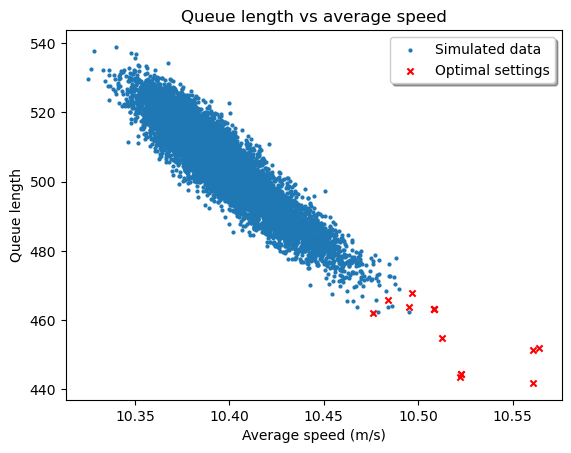}
\caption{\label{fig:distr_w_opt} Scatter plot with the average speed on the $x$-axis and the corresponding queue length on the $y$-axis. In blue are the results for the configurations used for training, while in red are the results for the optimal configurations according to the different surrogate models. It is clear that most of the optimal configurations are indeed better than the best configurations in the training set.}
\end{center}
\end{figure}

In Figure \ref{fig:distr_w_opt} we can see how the optimal configurations found using the Differential Evolution algorithm compare with the random configurations that were used to generate the training data. Blue dots represent traffic congestion variables from training configurations, while red dots correspond to the 14 optimal configurations described earlier (for which the results are presented in Tables \ref{tab:results_speed}-\ref{tab:results_both}). We can notice that the optimization procedure is indeed able to suggest configurations that significantly improve the traffic conditions.    

\section{Conclusions}\label{sec:conclusions}

This paper introduced a new macroscopic model which is an extension of the well-known Payne-Whitham model to the case of multiple intersections with traffic signals. The behaviour of cars in the presence of traffic lights was modelled by imposing suitable coupling conditions to the Payne-Whitham system of ODEs on the metric graph that describes the road structure. Another innovation we introduce involves incorporating the geometry of the road network to create more realistic coupling conditions that model the speed reduction occurring when cars turn at intersections. This speed reduction is dependent on the angle between the intersecting streets.

The resulting model (PWTL) was later used to find good configurations of traffic signals using the Differential Evolution optimization algorithm. However, in order to speed up the computations (which are time-consuming), the model was replaced by a variety of surrogate models, based on Linear Regression, Support Vector Regression, Multi Linear Perceptron and Random Forest. 

The SVR and RF models emerged as the most accurate surrogate models, achieving the lowest RMSE values. However, the Linear Regression model, despite being the least accurate in terms of RMSE on the training set, often yielded the best results in the simulated traffic dynamics, achieving the highest speeds and the lowest queue lengths. SVR, MLP and Random Forest, while accurate as surrogate models, did not consistently produce the best traffic configurations in simulations. These findings suggest that while more complex models are effective in replicating the PWTL framework, simpler models like Linear Regression can still provide valuable insights and optimal solutions in practical applications. In general, almost all surrogate models were able to suggest traffic light configurations that optimized speed and queue length, as highlighted in Figure \ref{fig:distr_w_opt}, demonstrating their effectiveness in improving traffic conditions. This ability to identify optimal configurations highlights the practical utility of surrogate models in reducing the computational burden associated with direct simulation of the PWTL framework. This is also consistent with similar research carried out using microscopic traffic simulations, surrogate models based on graph neural networks, and gradient descent as an optimization algorithm \cite{skowronek2021graph,skowronek2023graph}.

Moving forward, we aim to investigate deeper into the peculiarities of optimal traffic light configurations to better understand the factors that contribute to improved traffic flow. By examining these configurations in greater detail, we can identify specific patterns and strategies that lead to enhanced traffic management. Additionally, we plan to increase the training set to enhance the precision of our surrogate models, ensuring more accurate predictions and better optimization results.

Also, we will think about investigating other optimization techniques and comparing them with the Differential Evolution method.

Another important aspect is the observability of traffic characteristics, since it is usually not possible or very difficult to obtain all the necessary traffic data. Therefore, it is important to investigate how using only partial knowledge of the real traffic can affect the quality of the results.

\noindent\section*{Acknowledgments}

This article is based upon work from COST Action CA18232 MAT-DYN-NET, supported by COST (European Cooperation in Science and Technology). We would also like to thank the researchers from the Traffic Technical Institute at the Faculty of Civil and Geodetic Engineering, University of Ljubljana, especially Irena Strand, for providing and preparing the data that we used for calibrating the fundamental diagram of traffic flow.

D. Manea has been partially supported by the PNRR-III-C9-2023-I8 grant CF 149/31.07.2023 “Conformal Aspects of Geometry and Dynamics”.

\appendix
\section{Surrogate models parameters}\label{app:surrogate_par}

In our study, we used the Python \code{scikit-learn} package to evaluate multiple regression models with carefully chosen hyperparameters to identify the most effective surrogate model for our multi-output prediction task. We included four different models: Linear Regression, Support Vector Regression (SVR), Multi Linear Perceptron (MLP), and  Random Forest. For each of the models, we employed a 5-fold cross-validation approach and used the Halving Grid Search algorithm \cite{jamieson2016non} to find the optimal hyperparameters. 

For the Linear Regression model, we employed the default configuration provided by the \code{LinearRegression} class, which does not require hyperparameter tuning.

The \code{SVR} model was configured with an epsilon of 0.01 and a polynomial degree of 3, utilizing an extended cache size of 800 MB to accommodate larger datasets. The hyperparameters tuned were the kernel type (\code{rbf} and \code{sigmoid}), the gamma parameter (0.5 and 0.6), and the regularization parameter \code{C} (0.86 and 0.88). These parameters were selected to explore different kernel functions and their respective complexities, ensuring a thorough investigation of the model’s flexibility and generalization capabilities.

For the \code{MLPRegressor} model, we utilized a configuration with a hidden layer structure of (100, 50), examining both \code{tanh} and \code{relu} activation functions. The solver was fixed to \code{adam}, and we varied the alpha parameter (0.0001 and 0.001) for regularization. The learning rate was also tuned (0.001 and 0.01), along with a fixed batch size of 32. Early stopping was enabled to prevent overfitting, and the maximum number of iterations was set to 500. This setup allows the exploration of different network depths and activation dynamics, balancing model capacity and training efficiency.

Lastly, the \code{RandomForestRegressor} model was evaluated with different numbers of trees (50, 100, 200) and maximum depths (\code{None}, 5, 10, 20) to control the complexity of individual trees. The minimum samples required to split a node were set to 2 and 5. This range of parameters aims to optimize the trade-off between model variance and bias, ensuring robust performance across varied dataset complexities.

\printbibliography
\end{document}